\documentclass[10pt]{article}
\usepackage{amsmath}
\usepackage{latexsym}
\usepackage{amssymb}
\usepackage{amsfonts}
\usepackage{amsthm}
\title{THE CONFORMAL GROUP REVISITED}
\author{J.-F. POMMARET \\ CERMICS, Ecole des Ponts ParisTech\\ 
 jean-francois.pommaret@wanadoo.fr \\
http://cermics.enpc.fr/$\sim$pommaret/home.html }
\date{  }
\begin{document}
\maketitle

\noindent
{\bf ABSTRACT}  \\

\noindent
Since $100$ years or so, it has been usually accepted that the {\it conformal group} could be defined in an arbitrary dimension $n$ as the group of transformations preserving a non degenerate flat metric up to a nonzero invertible point depending factor called " {\it conformal factor} ". However, when $n\geq 3$, it is a finite dimensional Lie group of transformations with $n$ {\it translations}, $n(n-1)/2$ {\it rotations}, $1$ {\it dilatation} and $n$ nonlinear transformations called {\it elations} by E. Cartan in $1922$, that is a total of $(n+1)(n+2)/2$ transformations. Because of the Michelson-Morley experiment, the conformal group of space-time with $15$ parameters is well known for the Minkowski metric and is the biggest group of invariance of the Minkowski constitutive law of electromagnetism (EM) in vacuum, even though the two sets of field and induction Maxwell equations are respectively invariant by any local diffeomorphism. As this last generic number is also well defined and becomes equal to $3$ for $n=1$ or $6$ for $n=2$, the purpose of this paper is to use modern mathematical tools such as the {\it Spencer operator} on systems of OD or PD equations, both with its restriction to their symbols leading to the {\it Spencer} $\delta$-{\it cohomology}, in order to provide a {\it unique striking definition} that could be valid for any $n\geq 1$. The concept of an " {\it involutive system} " is crucial for such a new definition.   \\

\vspace{4cm}


\newpage

\noindent
{\bf 1) INTRODUCTION}   \\
Using local notations, this paper is mainly concerned with the following two connected problems:  \\
Given a differential operator $\xi \stackrel{{\cal{D}}}{\longrightarrow} \eta$, how can we find {\it compatibilty conditions} (CC), that is how can we construct a sequence $\xi \stackrel{{\cal{D}}}{\longrightarrow} \eta \stackrel{{\cal{D}}_1}{\longrightarrow} \zeta$ such that ${\cal{D}}_1 \circ {\cal{D}}=0$ and, among all such possible sequences, what are the " {\it best} "  ones, at least among the generating ones and when could we say that the sequence obtained is " {\it exact} " in a purely formal way, that is using only computer algebra for testing such a property. The order of an operator will be indicated under its arrow. \\
The difficulty is that, physicists being more familiar with analysis, they will say that a sequence is " {\it locally exact} " if one can find locally $\xi$ such that ${\cal{D}}\xi = \eta$ whenever ${\cal{D}}_1 \eta=0$. However, they have in mind the property of the exterior derivative $d$ and Maxwell equations in electromagnetism (EM), that is to say, using standard notations, the (local) possibility to introduce the EM  potential $A$ such that $dA=F$ whenever the EM field $F$ is a closed $2$-form with $dF=0$. \\
The main purpose of this paper is to prove that " things " could be much more delicate as these problems are only rarely associated with exterior calculus. We use the notations that can be found at length in our many books ([9],[10],[11],[12],[13],[14]) or papers ([15],[16],[18],[23],[24]).\\

 Let $\mu=({\mu}_1,...,{\mu}_n)$ be a multi-index with {\it length} ${\mid}\mu{\mid}={\mu}_1+...+{\mu}_n$, {\it class} $i$ if ${\mu}_1=...={\mu}_{i-1}=0,{\mu}_i\neq 0$ and $\mu +1_i=({\mu}_1,...,{\mu}_{i-1},{\mu}_i +1, {\mu}_{i+1},...,{\mu}_n)$. We set $y_q=\{y^k_{\mu}{\mid} 1\leq k\leq m, 0\leq {\mid}\mu{\mid}\leq q\}$ with $y^k_{\mu}=y^k$ when ${\mid}\mu{\mid}=0$. If $E$ is a vector bundle over $X$ with local coordinates $(x^i,y^k)$ for $i=1,...,n$ and $k=1,...,m$, we denote by $J_q(E)$ the $q$-{\it jet bundle} of $E$ with local coordinates simply denoted by $(x,y_q)$ and {\it sections} ${\xi}_q:(x)\rightarrow (x,{\xi}^k(x), {\xi}^k_i(x), {\xi}^k_{ij}(x),...)$ transforming like the section $j_q(\xi):(x)\rightarrow (x,{\xi}^k(x),{\partial}_i{\xi}^k(x),{\partial}_{ij}{\xi}^k(x),...)$ when $\xi$ is an arbitrary section of $E$. Then both ${\xi}_q\in J_q(E)$ and $j_q(\xi)\in J_q(E)$ are over $\xi \in E$ and the {\it Spencer operator}, {\it which is defined on sections}, just allows to distinguish them by introducing a kind of "{\it difference}" through the operator $d:J_{q+1}(E)\rightarrow T^*\otimes J_q(E): {\xi}_{q+1}\rightarrow j_1({\xi}_q)-{\xi}_{q+1}$ with local components $({\partial}_i{\xi}^k(x)-{\xi}^k_i(x), {\partial}_i{\xi}^k_j(x)-{\xi}^k_{ij}(x),...) $ and more generally $(d{\xi}_{q+1})^k_{\mu,i}(x)={\partial}_i{\xi}^k_{\mu}(x)-{\xi}^k_{\mu+1_i}(x)$. In a symbolic way, {\it when changes of coordinates are not involved}, it is sometimes useful to write down the components of $d$ in the form $d_i={\partial}_i-{\delta}_i$. The restriction of $d$ to the kernel $S_{q+1}T^*\otimes E$ of the canonical projection ${\pi}^{q+1}_q:J_{q+1}(E)\rightarrow J_q(E)$ is {\it minus} the {\it Spencer map} $\delta=dx^i\wedge {\delta}_i:S_{q+1}T^*\otimes E\rightarrow T^*\otimes S_qT^*\otimes E$ and $\delta \circ \delta =0$. The kernel of $d$ is made by sections such that ${\xi}_{q+1}=j_1({\xi}_q)=j_2({\xi}_{q-1})=...=j_{q+1}(\xi)$. Finally, if $R_q\subset J_q(E)$ is a {\it system} of order $q$ on $E$ locally defined by linear equations ${\Phi}^{\tau}(x,y_q)\equiv a^{\tau\mu}_k(x)y^k_{\mu}=0$, the $r$-{\it prolongation} $R_{q+r}={\rho}_r(R_q)=J_r(R_q)\cap J_{q+r}(E)\subset J_r(J_q(E))$ is locally defined when $r=1$ by the set of linear equations ${\Phi}^{\tau}(x,y_q)=0, d_i{\Phi}^{\tau}(x,y_{q+1})\equiv a^{\tau\mu}_k(x)y^k_{\mu+1_i}+{\partial}_ia^{\tau\mu}_k(x)y^k_{\mu}=0$ and has {\it symbol} $g_{q+r}=R_{q+r}\cap S_{q+r}T^*\otimes E\subset J_{q+r}(E)$ if one looks at the {\it top order terms}. If ${\xi}_{q+1}\in R_{q+1}$ is over ${\xi}_q\in R_q$, differentiating the identity $a^{\tau\mu}_k(x){\xi}^k_{\mu}(x)\equiv 0$ with respect to $x^i$ and substracting the identity $a^{\tau\mu}_k(x){\xi}^k_{\mu+1_i}(x)+{\partial}_ia^{\tau\mu}_k(x){\xi}^k_{\mu}(x)\equiv 0$, we obtain the identity $a^{\tau\mu}_k(x)({\partial}_i{\xi}^k_{\mu}(x)-{\xi}^k_{\mu+1_i}(x))\equiv 0$ and thus the restriction $d:R_{q+1}\rightarrow T^*\otimes R_q$ ([9],[11],[12],[27]).   \\

\noindent
{\bf DEFINITION 1.1}: $g_q$ is said to be $s$-{\it acyclic} if the purely algebraic $\delta$-cohomology $H^s_{q+r}(g_q)$ of $... \stackrel{\delta}{\rightarrow} {\wedge}^sT^*\otimes g_{q+r}\stackrel{\delta}{\rightarrow}...$ are such that $H^1_{q+r}(g_q)= ... = H^s_{q+r}(g_q)=0,  \forall r\geq 0$ and {\it involutive} if it is $n$-acyclic. Also $R_q$ is said to be {\it involutive} if it is {\it formally integrable} (FI), that is when the restriction ${\pi}^{q+1}_q:R_{q+1}\rightarrow R_q $ is an epimorphism $\forall r\geq 0$ or, equivalently, when all the equations of order $q+r$ are obtained by $r$ prolongations only, $\forall r\geq 0$ and $g_q$ is involutive. In that case, $R_{q+1}\subset J_1(R_q)$ is a canonical equivalent formally integrable first order system on $R_q$ with no zero order equations, called the {\it Spencer form}.   \\

\noindent
{\bf EXAMPLE 1.2}: ({\it Classical Killing operator})  \\
Considering the {\it classical Killing} operator ${\cal{D}}:\xi \rightarrow {\cal{L}}(\xi)\omega=\Omega \in S_2T^*=F_0$ where ${\cal{L}}(\xi)$ is the Lie derivative with respect to $\xi$ and $\omega \in S_2T^*$ is a nondegenerate metric with $det(\omega)\neq 0$. Accordingly, it is a lie operator with ${\cal{D}}\xi=0, {\cal{D}}\eta=0 \Rightarrow {\cal{D}}[\xi,\eta]=0$ and we denote simply by $\Theta \subset T$ the set of solutions with $[\Theta,\Theta ]\subset \Theta$. Now, as we have explained many times, the main problem is to describe the CC of ${\cal{D}}\xi=\Omega \in F_0$ in the form ${\cal{D}}_1\Omega=0$ by introducing the so-called {\it Riemann} operator ${\cal{D}}_1:F_0 \rightarrow F_1$. We advise the reader to follow closely the next lines and to imagine why it will not be possible to repeat them for studying the {\it conformal Killing} operator. Introducing the well known Levi-Civita isomorphism $j_1(\omega)=(\omega, {\partial}_x \omega) \simeq (\omega , \gamma)$ by defining the Christoffel symbols 
${\gamma}^k_{ij}=\frac{1}{2}{\omega}^{kr}({\partial}_i{\omega}_{rj} + {\partial}_j{\omega}_{ir}-
{\partial}_r{\omega}_{ij})$ where $({\omega}^{rs})$ is the inverse matrix of $({\omega}_{ij})$ and the {\it formal Lie derivative}, we get the second order system $R_2 \subset J_2(T)$:  \\
\[  \left\{  \begin{array}{lcccl}
{\Omega}_{ij}& \equiv &(L({\xi}_1)\omega)_{ij}  & = & {\omega}_{rj}(x){\xi}^r_i + {\omega}_{ir}(x){\xi}^r_j + {\xi}^r {\partial}_r {\omega}_{ij}(x) =0  \\
{\Gamma}^k_{ij}& \equiv & (L({\xi}_2)\gamma)^k_{ij} & =  & {\xi}^k_{ij} + {\gamma}^k_{rj}(x) {\xi}^r_i+{\gamma}^k_{ir}(x) {\xi}^r_j +{\gamma}^k_{ir}(x) {\xi}^r_j - {\gamma}^r_{ij}(x) {\xi}^k_r + {\xi}^r{\partial}_r {\gamma }^k_{ij}(x)=0
\end{array} \right.  \]
with sections ${\xi}_2:x \rightarrow ({\xi}^k(x), {\xi}^k_i(x), {\xi}^k_{ij}(x))$ transforming like $j_2(\xi): x \rightarrow ({\xi}^k(x), {\partial}_i{\xi}^k(x), {\partial}_{ij}{\xi}^k(x))$. The system $R_1 \subset J_1(T)$ has a symbol $g_1 \simeq {\wedge}^2T^*\subset T^* \otimes T$ depending only on $\omega$ with $dim(g_1)=n(n-1)/2$ and is finite type because its first prolongation is $g_2=0$. It cannot be thus involutive and we need to use one additional prolongation. Indeed, using one of the main results to be found in ([12],[19],[23]), we know that, when $R_1$ is FI, then the CC of ${\cal{D}}$ are of order $s+1$ where $s$ is the number of prolongations needed in order to get a $2$-acyclic symbol, that is $s=1$ in the present situation, a result that should lead to CC of order $2$ if $R_1$ were FI. However, it is known that $R_2$ is FI, thus involutive, {\it if and only if} $\omega$ has constant Riemannian curvature, a result first found by L.P. Eisenhart in 1926 ([2]) which is only a particular example of the {\it Vessiot structure equations} discovered b E. Vessiot in $1904$ ([28]), though in a quite different setting (See [12] and [19] for an explicit modern proof). Such a necessary condition for constructing an exact differential sequence could not have been used by any follower because the " {\it Spencer machinery} " has only been known after 1970 ([27]). Otherwise, if the metric does not satisfy this condition, CC may exist but have no link with the Riemann tensor ([23]). We may define the vector bundle $F_1$ in the short exact sequence made by the top row of the following commutative diagram: \\
\[  \begin{array}{rcccccccl}
 &  &  & 0  & &  0 &  &  &  \\
 &  &  & \downarrow &  & \downarrow &  & &   \\
 &0 & \rightarrow &S_3T^*\otimes T &\rightarrow &S_2T^*\otimes F_0& \rightarrow & F_1 &  \rightarrow 0 \\
&  &  & \downarrow &  & \downarrow & & &  \\
 & 0 & \rightarrow &T^* \otimes S_2T^*\otimes T &\rightarrow &T^*\otimes T^*\otimes F_0& \rightarrow & 0 &  \\
& \downarrow  &  & \downarrow &  & \downarrow &  &  &  \\
    0\rightarrow &{\wedge}^2T^*\otimes g_1  & \rightarrow &\underline{{\wedge}^2T^* \otimes T^*\otimes T }  &\rightarrow &{\wedge}^2T^*\otimes  F_0& \rightarrow &0 &    \\
    &\downarrow   &  & \downarrow &  & \downarrow  &  &   \\
0 \rightarrow & {\wedge}^3T^* \otimes T & = & {\wedge}^3T^*\otimes T  & \rightarrow & 0& & &   \\
  &\downarrow   &  & \downarrow &  & &  & &   \\
   & 0  &  & 0 &  & &  & &   \\
 \end{array}   \]
where the vertical $\delta$-sequences are exact {\it but the first}, or, using a snake type diagonal chase, from the short exact sequence of vector bundles: \\
\[        0 \rightarrow F_1  \rightarrow T^*\otimes g_1 \stackrel{\delta}{\longrightarrow} {\wedge}^2T^*\otimes T \rightarrow 0  \]
This result is first leading to the long exact sequence of vector bundles:  \\
\[  0 \rightarrow     R_3   \rightarrow      J_3(T)    \rightarrow J_2(F_0) \rightarrow F_1  \rightarrow 0  \]
and to the {\it Riemann} operator ${\cal{D}}_1:F_0 \stackrel{j_2}{\longrightarrow} J_2(F_0) \rightarrow F_1$. As $g_2=0$, we also discover that $F_1$ is just the Spencer $\delta$-cohomology $H^2(g_1)$ at ${\wedge}^2T^*\otimes g_1$ along the previous short exact sequence.  \\
We get the {\it striking formulas} where the $+$ signs are replaced by $-$ signs:  \\
\[   \begin{array}{ rcl}
dim(F_1) &  =  &  n^2(n+1)^2/4 - n^2(n+1)(n+2)/6             \\
                 &  =  &  n^2(n-1)^2/4 -n^2(n-1)(n-2)/6     \\
                 &  =   &   n^2(n^2-1)/12 
  \end{array}       \]
This result, first found as early as in 1978 ([9]), clearly exhibit {\it without indices} the two well known algebraic properties of the Riemann tensor as a section of the tensor bundle ${\wedge}^2T^*\otimes T^* \otimes T$. \\

It thus remains to exhibit the {\it Bianchi} operator exactly as we did for the {\it Riemann} operator, with the same historical comments already provided. However, now we know that $R_1$ is formally integrable (otherwise nothing could be achieved and we should start with a smaller system [23]), the construction of the linearized Janet-type differential sequence as a strictly exact differential sequence but {\it not} an involutive differential sequence because the system $R_1$ and thus the first order operator ${\cal{D}}$ are formally integrable though {\it not} involutive as $g_1$ is finite type with $g_2=0$ but not involutive. Doing one more prolongation only, we obtain the first order {\it Bianchi} operator 
${\cal{D}}_2: F_1 \stackrel{j_1}{\longrightarrow} J_2(F_1) \rightarrow F_2$ as before, defining the vector bundle $F_2$ in the long exact sequence made by the top row of the following commutative diagram:  \\
\[  \begin{array}{rccccccccccl}
 &  &  & 0  & &  0 &  & 0 & & & \\
 &  &  & \downarrow &  & \downarrow &  & \downarrow & & &  \\
 &0 & \rightarrow &S_4T^*\otimes T &\rightarrow &S_3T^*\otimes F_0& \rightarrow &T^* \otimes F_1 & \rightarrow &F_2 & \rightarrow 0    \\
  &  &  & \downarrow &  & \downarrow &  & \parallel & & &  \\
 & 0 & \rightarrow &T^* \otimes S_3T^*\otimes T &\rightarrow &T^*\otimes S_2T^*\otimes F_0& \rightarrow &T^* \otimes F_1 & \rightarrow &0&   \\
  &  &  & \downarrow &  & \downarrow &  & \downarrow & & &  \\
 & 0 & \rightarrow &{\wedge}^2T^* \otimes S_2T^*\otimes T &\rightarrow &{\wedge}^2T^*\otimes T^*\otimes F_0& \rightarrow &0 & &&   \\
  &\downarrow   &  & \downarrow &  & \downarrow &  & & & &  \\
   0\rightarrow &{\wedge}^3T^*\otimes g_1  & \rightarrow &\underline{{\wedge}^3T^* \otimes T^*\otimes T} &\rightarrow &{\wedge}^3T^*\otimes  F_0& \rightarrow &0 & &&   \\
    &\downarrow   &  & \downarrow &  & \downarrow  &  & & & &  \\
0 \rightarrow & {\wedge}^4T^* \otimes T & = & {\wedge}^4T^*\otimes T  & \rightarrow & 0& & & & &  \\
  &\downarrow   &  & \downarrow &  & &  & & & &  \\
   & 0  &  & 0 &  & &  & & & &  \\
 \end{array}   \]
where the vertical $\delta$-sequences are exact {\it but the first}, or, using a snake type diagonal chase, from the short exact sequence:  \\
\[   0 \rightarrow F_2 \rightarrow {\wedge}^3T^*\otimes g_1 \stackrel{\delta}{\longrightarrow} {\wedge}^4T^*\otimes T  \rightarrow 0  \]
showing that $F_2=H^3(g_1) $ ([8],[9]). We have in particular for $n\geq 4$:   \\
\[ \begin{array}{ lcl}
dim(F_2) & = & n^2(n-1)^2(n-2) /12 - n^2(n-1)n-2)(n-3)/24  \\
     & = &  n^2(n+1)(n+2(n+3)/24 +n^3(n^2_1)/12 - n^2(n+1)(n+2)(n+3)/24  \\
& = & n^2(n^2 - 1)(n-2)/24
\end{array}    \]
and thus $dim(F_2)=(4\times 6)-(1\times 4)=(16\times 15 \times 2)/24=20$ when $n=4$. This result also exhibit all the properties of the Bianchi identities as a section of the tensor bundle 
${\wedge}^3 T^* \otimes T^* \otimes T$. In arbitrary dimension, we finally obtain the differential sequence, which is {\it not} a Janet sequence:  \\
\[  0 \rightarrow \Theta \rightarrow  T  \underset 1{\stackrel{Killing}{\longrightarrow}}F_0  \underset 2{ \stackrel{Riemann}{\longrightarrow}} F_1  \underset 1{\stackrel{Bianchi}{\longrightarrow}} F_2    \]  

\noindent
{\bf EXAMPLE 1.3}: ({\it Conformal Killing operator})  \\
At first sight, it seems that similar methods could work in order to study the conformal Killing operator 
and, more generally, all conformal concepts will be described with a "hat", in order to provide the strictly exact differential sequence:  \\
\[  0 \rightarrow \hat{\Theta} \rightarrow  T  \stackrel{\hat{\cal{D}}}{\longrightarrow} {\hat{F}}_0   \stackrel{{\hat{\cal{D}}}_1}{\longrightarrow} 
{\hat{F}}_1  \stackrel{{\hat{\cal{D}}}_2}{\longrightarrow} {\hat{F}}_2    \]  
where ${\hat{\cal{D}}}_1$ is the $Weyl$ operator with generating CC ${\hat{\cal{D}}_2}$. It is only in $2016$ (See [18] for more details) that we have been able to recover all these operators and confirm with computer algebra that the orders of the operators involved highly depend on the dimension, even when $n\geq 3$ as follows:  \\

\noindent 
 $\bullet$ $n=3$: $ \hspace{2cm} 3 \underset{1}{\longrightarrow} 5 \underset{3}{\longrightarrow }5 \underset{1}{ \longrightarrow }3 \rightarrow 0 $ \\   \\ 
 $\bullet$ $n=4$: $ \hspace{2cm} 4 \underset{1}{\longrightarrow} 9 \underset{2}{\longrightarrow} 10 \underset{2}{\longrightarrow} 9 \underset{1}{\longrightarrow} 4  \rightarrow 0  $ \\    \\                                       
 $\bullet$ $n\geq 5$: $ \hspace{2cm}  5 \underset{1}{\longrightarrow} 14 \underset{2}{\longrightarrow} 35 \underset{1}{\longrightarrow} 35 \underset{2}{\longrightarrow} 14 \underset{1}{\longrightarrow} 5 \rightarrow 0  $\\

This result is based on the following technical lemma (See [19] for details):  \\

\noindent
{\bf LEMMA 1.4}: The symbol ${\hat{g}}_1$ defined by the linear equations:  \\
\[   {\hat{\Omega}}_{ij}\equiv {\omega}_{rj}(x) {\xi}^r_i + {\omega}_{ir}(x){\xi}^r_j - \frac{1}{2}{\omega}_{ij}(x){\xi}^r_r=0  \]
does not depend on any conformal factor, is finite type with ${\hat{g}}_3=0, \forall n\geq 3$ and is surprisingly such that ${\hat{g}}_2$ is $2$-acyclic 
for  $n\geq 4$ or even $3$-acyclic when $n\geq 5$.   \\

\noindent
{\bf REMARK 1.5}: In order to emphasize the reason for using Lie equations, we now provide the explicit form of the $n$ infinitesimal elations with $ 1\leq r,s,t \leq n$, whenever $n\geq 3$:   \\
\[  {\theta}_ s = - \frac{1}{2} x^2 {\delta}^r_s{\partial}_r+{\omega}_{st}x^tx^r{\partial}_r   \, \Rightarrow  \, {\partial}_r{\theta}^r_s=n{\omega}_{st}x^t, \,\, [{\theta}_s,{\theta}_t]=0  \]
where the underlying metric is used for the scalar product $x^2$ involved. It is easy to check that ${\xi}_2 \in S_2T^*\otimes T$ defined by $  {\xi}^k_{ij}(x)= {\lambda}^s(x){\partial}_{ij} {\theta}^k_s(x)$ belongs to ${\hat{g}}_2$ with $A_i={\omega}_{si}{\lambda}^s$ in the following formula where $\delta$ is the standard Kronecker symbol and ${\xi}_2 \in {\hat{R}}_2$:  \\
\[   {\Gamma}^k_{ij} \equiv (L({\xi}_2)\gamma)^k_{ij}= {\xi}^k_{ij} + {\gamma}^k_{rj}{\xi}^r_i + {\gamma}^k_{ir}{\xi}^r_j -  {\gamma}^r_{ij} {\xi}^k_r + {\xi}^r {\partial}_r {\gamma}^k_{ij} = {\delta}^k_iA_j + {\delta}^k_jA_i - {\omega}_{ij}{\omega}^{kr} A_r  \] 
We thus understand how important it is to use " {\it sections} "  rather than " {\it solutions} ".  \\

Accordingly, a possible unification can be achieved through the " {\it fundamental diagram I} " relating together the {\it Spencer sequence} and the {\it Janet sequence} as follows in arbitrary dimension $n$ for any {\it involutive} system $R_q\subseteq J_q(E)$ because these are the only existing canonical sequences:  \\
 
  \[  \begin{array}{rccccccccccccccl}
 &&&&& 0 &&0&&0&  & 0 & & \\
 &&&&& \downarrow && \downarrow && \downarrow  & & \downarrow  &  &\\
  & 0& \rightarrow& \Theta &\stackrel{j_q}{\rightarrow}& C_0  &\stackrel{D_1}{\rightarrow}&  C_1 &\stackrel{D_2}{\rightarrow} ... \stackrel{D_{n-1}}{\rightarrow}&  C_{n-1} &\stackrel{D_n}{\rightarrow}  & C_n & \rightarrow  & 0  \\
  &&&&& \downarrow & & \downarrow & & \downarrow & & \downarrow &   & \\
   & 0 & \rightarrow &  E & \stackrel{j_q}{\rightarrow} & C_0(E)  & \stackrel{D_1}{\rightarrow} &  C_1(E)  &\stackrel{D_2}{\rightarrow}... \stackrel{D_{n-1}}{\rightarrow} & C_{n-1}(E) &  \stackrel{D_n}{\rightarrow} & C_n(E) & \rightarrow  & 0 \\
   & & & \parallel && \hspace{5mm}\downarrow {\Phi}_0 & &\hspace{5mm} \downarrow {\Phi}_1 & & \hspace{5mm}\downarrow {\Phi}_{n-1} &  &\hspace{5mm}\downarrow {\Phi}_n & \\
   0 \rightarrow & \Theta &\rightarrow &  E & \stackrel{\cal{D}}{\rightarrow} &  F_0  & \stackrel{{\cal{D}}_1}{\rightarrow} &  F_1  & \stackrel{{\cal{D}}_2}{\rightarrow} ... \stackrel{{\cal{D}}_{n-1}}{\rightarrow}&  F_{n-1} & \stackrel{{\cal{D}}_n}{\rightarrow}&  F_n  & \rightarrow  &  0      \\
   &&&&& \downarrow & & \downarrow & & \downarrow &  & \downarrow &    &   \\
   &&&&& 0 && 0 && 0  &  &  0  &  &
   \end{array}     \]
where  $C_0=R_q \subset J_q(E)=C_0(E)$ and $dim(F_r)= dim(C_r(E)) - dim (C_r)$. Indeed, we have $dim(C_r)= dim({\wedge}^rT^*)\times dim(R_q)$ for finite type involutive systems and we therefore  notice that {\it the crucial point is to deal with involutive systems}. In the group framework, we have $E=T$ and, as we are dealing with finite type systems, it is thus sufficient to replace $j_q$ and $R_q\subset J_qE)$ by $j_2$ and $R_2\subset J_2(T)$ with $g_2=0$ in the classical situation or by $j_3$ and 
${\hat{R}}_3 \subset J_3(T)$ with ${\hat{g}}_3=0$ in the conformal situation, on the ondition to be able to treat the specific cases $n=1$ and $n=2$.  \\  

Finally, as a different way to look at these questions, if $K$ be a differential field containing $\mathbb{Q}$, we may introduce the ring $D=K[d]=K[d_1,...,d_n]$ of differential operators with coefficients in $K$ and consider a linear differential operator ${\cal{D}}$ with coefficients in $K$. If ${\cal{D}}_1$ generates the CC of ${\cal{D}}$, we have of course ${\cal{D}}_1 \circ {\cal{D}} = 0$. Taking the respective (formal) {\it adjoint} operators, we obtain therefore $ad({\cal{D}}) \circ ad({\cal{D}}_1)=0$ but $ad({\cal{D}})$ may not generate the CC of $ad({\cal{D}}_1)$ and so on in {\it any} differential sequence where each operator generates the CC of the preceding one.  \\

\noindent
{\bf DEFINITION 1.6}: If $M$ is the differential module over $D$ or simply $D$-module defined by ${\cal{D}}$, we set $ext^0_D(M)=hom_D(M,D)$. As for the other {\it extension modules}, they have been created in order to " {\it measure} " the previous gaps ([7],[12],[24]). In particular, we say that $ext^1_D(M)=0$ if  $ad({\cal{D}})$ generates the CC of $ad({\cal{D}}_1)$, that  $ext^2_D(M)=0$ if $ad({\cal{D}}_1)$ generates the CC of $ad({\cal{D}}_2)$ and so on. Moreover, if ${\cal{D}}$ is of finite type, then $ad({\cal{D}})$ is surjective with $ext^0_D(M)=0$. The simplest example is that of classical space geometry with $n=3$ and $ad(grad)= - div$. Similar definitions are also valid for the Janet and Spencer sequences. Also, vanishing of the first extension module amounts to the existence of a local parametrization by potential-like functions. \\
 
According to a (difficult) theorem of (differential) homological algebra, {\it the extension modules only depend on $M$ and not on the previous differential sequences used} ([8],[25]). They are used in agebraic geometry and have even been introduced in engineering sciences after $1990$ (control theory) ([13],[14],[24]). However, though the extension modules are the {\it only} intrinsic objects that can be associated with a differential module, they have surprisingly {\it never} been introduced in mathematical physics. The main problem is that a control system is controllable if and only if it is parametrizable by potentials while the systems involved can be parametrized in all classical physics (Cauchy or Maxwell equations are well known examples in [15]) apart from ... Einstein equations ([16],[18],[22],[24]). As for the tools involved, we let the reader compare ([3],[6]) to ([10],[11]).  \\

After presenting two motivating examples in Section $2$, such a procedure will be achieved in Section 
$3$ in such a way that the Spencer sequences involved, being isomorphic to tensor products of the Poincar\'{e} sequence by finite dimensional Lie algebras, will have therefore vanishing zero, first and second extension modules when $n\geq 3$ ([12],[24]). For all results concerning differential modules, we refer the reader to the (difficult) references ([1],[5],[26]) or to ([12],[22]).  \\    \\

\noindent
{\bf 2) TWO MOTIVATING EXAMPLES}  \\
\noindent
{\bf EXAMPLE 2.1}  \\
With $m=1, n=2, q=2, K= \mathbb{Q}$, let us consider the inhomogeneous second order operator:   \\
\[      Py\equiv d_{22}y=u , \hspace{3cm}    Qy\equiv d_{12}y-y=v   \]
We obtain at once through crossed derivatives:\\
\[      y=d_{11}u-d_{12}v-v  \,\,\,  \Rightarrow  \,\,\,  \Theta = 0    \]
and, by substituting, two fourth order CC for $(u,v)$, namely:\\
\[   \left \{   \begin{array}{lclcl}
U & \equiv & d_{1122}u-d_{1222}v-d_{22}v-u &=& 0\\
V & \equiv & d_{1112}u-d_{11}u-d_{1122}v & = & 0
\end{array} 
\right \} \Rightarrow  W \equiv d_{12}V+V-d_{11}U=0 \] 
However, the commutation relation $P\circ Q\equiv Q\circ P$ provides a single CC for $(u,v)$, namely:\\
\[     C \equiv d_{12}u-u-d_{22}v=0   \]
and we check at once $   U=d_{12}C+C, V=d_{11}C  $  while $ C=d_{22}V-d_{12}U+U $, that is:\\
 \[(U=0, V=0) \Leftrightarrow (C=0).\]
 Using corresponding notations, let us compare the two following differential sequences:  \\
 \[  \begin{array}{lcc}
   0 \rightarrow \Theta \rightarrow y \underset 2{\stackrel{{\cal{D}}}{\longrightarrow}} (u,v)\underset 4 {\stackrel{{\cal{D}}_1}{\longrightarrow}} (U,V) \underset 2 {\stackrel{{\cal{D}}_2}{\longrightarrow}} W  \rightarrow 0  &\hspace{2cm} &  (1)  \\
 0 \rightarrow \Theta \rightarrow y \underset 2 {\stackrel{{\cal{D}}}{\longrightarrow}} (u,v) \underset 2 {\stackrel{{\cal{D}}'_1}{\longrightarrow}} C \rightarrow  0  & \hspace{2cm}  &  (2)   
\end{array}   \]
Though the second order system considered is surely {\it not} FI because the $4$ parametric jets of $R_2$ are $(y, y_1, y_2, y_{11})$ and the $4$ (again !) parametric jets of $R_3$ are $(y, y_1, y_{11}, y_{111})$ but the $4$ (again !) parametric jets of $R_4$ are $(y_1, y_{11}, y_{111}, y_{1111})$. More generally, we let the reader prove by induction that $dim(R_{2+r})=4, \forall r\geq 0$. The formal $r$-prolongation of $(2)$, namely:   \\
\[    0 \rightarrow R_{r+4} \rightarrow J_{r+4}(y) \rightarrow J_{r+2}(u,v) \rightarrow J_r(C) \rightarrow 0  \]
is exact because $ 4 - (r+5)(r+6)/2 + (r+3)(r+4) - (r+1)(r+2)/2= 0$, even though the corresponding symbol sequence: \\
\[    0 \rightarrow g_{r+4} \rightarrow S_{r+4}T^*Ê(y) \rightarrow S_{r+2} T^*(u,v) \rightarrow S_rT^*(C) \rightarrow 0  \]
is {\it not} exact because $ ((2(r+3)) - (r+1)) - ((r+5) - 1) = (r+5) -(r+4) = 1\neq 0$ because the system considered is not formally integrable.  \\
On the contrary, the prolongations of $(1)$ are {\it not} exact on the jet level. Indeed, the long sequence:  \\
\[       0 \rightarrow R_8 \rightarrow J_8(y) \rightarrow J_6(u,v) \rightarrow J_2(U,V) \rightarrow W
 \rightarrow 0   \]
is {\it not} exact because we have $  4 - 45 + 56 - 12 +1 = 4 \neq 0  $. \\
Now, considering the ring $D=\mathbb{Q} [d_1,d_2]$ of differential operators with coefficients in the trivial differential field $\mathbb{Q}$, we have the "exact" sequences of differential modules where $M=0$:  \\
\[   \begin{array}{rcc}
  0\rightarrow D  \rightarrow D^2 \rightarrow D^2 \rightarrow D \stackrel{p}{\longrightarrow} M \rightarrow 0  & \hspace{2cm} &  (1^*)   \\
  0\rightarrow D  \rightarrow D^2 \rightarrow  D \stackrel{p}{\longrightarrow} M \rightarrow 0  & \hspace{2cm} &  (2^*)   \\
\end{array}   \]
where $ p $ is the canonical residual projection. However, and this is a quite delicate point rarely known even by mathematicians, {\it a fortiori} by physicists, they are not " {\it strictly} " exact even if the Euler-Poincar\'{e} characteristics both vanish because $1 - 2 + 2 - 1=0$ and $1  2 + 1 = 0$. (See [15] for definitions and more details). Roughly speaking, it follows that the " {\it best} " differential sequences are obtained by using {\it only} formally integrable operators/systems in such a way that sequences on the jet level can be studied through their symbol sequences, the " {\t canonical} " ones by using {\it exclusively} involutive operators/systems in such a way that what happens with ${\cal{D}}$ also hapens with ${\cal{D}}_1$ and so on. It follows that the sequences $(2)$ or $(2^*)$ are "{\it better} " than $(1)$ or $(1^*)$ because they provide {\it more informations} on the generating CC. \\
However, the given system is {\it not} FI and it should be " {\it better} " to use another system providing {\it more informations}. In particular, if we start wth a system $R_q\subset J_q(E)$ and set $R_{q+r}={\rho}_r(R_q)=J_r(R_q) \cap J_{q+r}(E)$, it is known that (in general) one can find two integers $r,s\geq 0$ such that the system $R^{(s)}_{q+r}={\pi}^{q+r+s}_{q+r}(R_{q+r+})$ is formally integrable and even involutive with the same solutions ([9],[13],[14]). When all the operators are FI, the sequence is said to be {\it strictly} exact ([17]).   \\
In the present situation, it should be " {\it better} " to replace $R_2$ by $R^{(4)}_2=0$ because $R^{(2)}_2$ is adding $y=0$ while $R^{(3)}_2$ is adding $y_1=0, y_2=0$ and $R^{(4)}_2$ is adding $y_{11}=0$. It follows that the Janet sequence for the injective trivially involutive operator $j_2$ is providing even {\it more informations}, along with the fact that the Spencer bundles vanish in the " {\it fundamental diagram I} " ([9],[12],[13]).  \\
We let the reader check that all the extension modules vanish because $M=0$ and to compare with the {\it totally different} involutive system defined by $y_{22}=0, y_{12}=0$ with $M\neq 0 \Rightarrow ext^0(M)\neq 0, ext^1(M)\neq  0, ext^2(M) \neq 0$.  \\

\noindent
{\bf  EXAMPLE 2.2}  \\
$\bullet$ {\it FIRST STEP} \,\, With $n=3, m=1, q=2$, let us consider the second order linear system $R_2 \subset J_2(E)$ introduced by F.S. Macaulay in his $1916$ book ([7])(See also [23] for more details):  \\
\[  {\Phi}^3\equiv y_{33}=0, \,\, {\Phi}^2\equiv y_{23} - y_{11}=0, \,\, {\Phi}^1\equiv y_{22}=0  \]
Using muli-indices, we may introduce the operators $R=d_{33}, \,\, Q=d_{23} - d_{11}, \,\, P=d_{22} $. Taking into account the $3$ commutation relations $[Q,R]=0, [R,P]=0, [P,Q]=0$ and the single Jacobi identity $[P,[Q,R]] + [Q,[R,P]] + [R,[P,Q]]=0, \forall (P,Q,R)$, we obtain at once the following locally and strictly exact sequence where the order of each operator is under its own arrow:   \\
\[  0 \longrightarrow \Theta \longrightarrow 1 \underset 2{\stackrel{{\cal{D}}}{\longrightarrow}} 3 \underset 2{\stackrel{{\cal{D}}_1}{\longrightarrow}} 3 \underset 2{\stackrel{{\cal{D}}_2}{\longrightarrow}} 1 \longrightarrow 0  \]
However, the first operator ${\cal{D}}$ involved cannot be involutive because it is finite type, that is $g_{q+r}=0$ for a certain integer $r\geq 0$ as we must have an exact sequence $0 \rightarrow {\wedge}^{(n-1)}T^*\otimes g_{q+r-1} \rightarrow 0$ and so on. The first prolongation is obtained by adding the $9$ PD equations:   \\
\[  y_{333}=0, y_{233}=0, y_{223}=0, y_{222}=0, y_{133}=0, y_{123} - y_{111}=0, y_{122}=0,
 y_{113}=0,y_{112}=0 \]
and the second prolongation is obtained by adding the $15$ PD equations $y_{ijkl}=0$. We obtain therefore $dim(g_2)=6 - 3=3, dim(g_3)=1, g_4=0$. Nevertheless, the interesting fact is that $g_3$ is 
$2$-acyclic without being $3$-acycic and thus involutive. Indeed, we have the $\delta$-sequences:  \\
\[   0 \rightarrow {\wedge}^2T^*\otimes g_3 \stackrel{\delta}{\longrightarrow} {\wedge}^3T^*\otimes g_2 \rightarrow 0 , \hspace{1cm}  0 \rightarrow {\wedge}^3T^*\otimes g_3 \rightarrow 0  \]
Using the letter $v$ for the symbol coordinates, the mapping $\delta$ on the left is defined by:   \\
\[      v_{111,23} + v_{112,31} + v_{113,12} =v_{11,123},  
         v_{121,23} + v_{122,31} + v_{123,12} =v_{12,123},                  
         v_{131,23} + v_{132,31} + v_{133,12} =v_{13,123}    \]
that is to say $v_{111,23}=v_{11,23}, \,\,\, v_{111,12}=v_{12,123}, \,\,\, v_{111,31}=v_{13,123}$. The corresponding $ \delta$-map is thus injective {\it and} surjective, that is $g_3$ is $2$-acyclic but cannot be also $3$-acyclic because of the inequality, $dim({\wedge}^3T^*\otimes g_3)=dim(g_3)=1\neq 0$. The above sequence is thus very far from being a Janet sequence and we cannot compare it with the Spencer sequence. \\

\noindent
$\bullet$ {\it SECOND STEP} \,\, In the example of Macaulay, we have at once $dim(R_2)=7$ with the $7$ parametric jets $(y, y_1, y_2, y_3, y_{11}, y_{12}, y_{13})$ and thus $dim(R_4)=dim(R_3)=7+1=8=2^3$ with the only additional third order parametric jet $(y_{111})$. We notice that, when $n=2$, the new system $R_2$ defined by $y_{22}=0, y_{12} - y_{11}=0$ is also finite type with $y_{ijr}=0$ and thus $dim(R_3)=dim(R_2)=4=2^2$ (See [7] for this striking result on the powers of $2$). Therefore, instead of starting with the previous second order operator ${\cal{D}}_1$ defined by $R_2$, we may now start afresh with the new third order operator ${\cal{D}}_1$ defined by $R_3$ which is not involutive again. We let the reader check as a tricky exercise or using computer algebra that one may obtain " {\it necessarily} " the following finite length differential sequence which is far from being a Janet sequence but for other reasons.  \\

\[ 0 \rightarrow \Theta \rightarrow E \underset 3{\stackrel{{\cal{D}}}{\longrightarrow}} F_0
                          \underset 1{\stackrel{{\cal{D}}_1}{\longrightarrow}} F_1
                         \underset 2{\stackrel{{\cal{D}}_2}{\longrightarrow}} F_2
                           \underset 1{\stackrel{{\cal{D}}_3}{\longrightarrow}} F_3
                        \underset 1{\stackrel{{\cal{D}}_4}{\longrightarrow}} F_4
                          \underset 1{\stackrel{{\cal{D}}_5}{\longrightarrow}}F_5 \rightarrow 0    \]

\[ 0 \rightarrow \Theta \rightarrow 1 \underset 3{\stackrel{{\cal{D}}}{\longrightarrow}} 12
                          \underset 1{\stackrel{{\cal{D}}_1}{\longrightarrow}} 21
                         \underset 2{\stackrel{{\cal{D}}_2}{\longrightarrow}} 46
                           \underset 1{\stackrel{{\cal{D}}_3}{\longrightarrow}} 72
                        \underset 1{\stackrel{{\cal{D}}_4}{\longrightarrow}} 48
                          \underset 1{\stackrel{{\cal{D}}_5}{\longrightarrow}}12 \rightarrow 0    \]                          
and we check that $1 - 12 + 21 - 46 + 72 - 48 + 12 = 0 $. As $g_3$ is $2$-acyclic, the third order operator 
${\cal{D}}$ has a CC operator ${\cal{D}}_1$ of order $1$ having a CC operator ${\cal{D}}_2$ of order 
$2$ which is involutive, {\it totally by chance}, and we end with the Janet sequence for ${\cal{D}}_2$. Such a situation is the only one we have met during the last ... $40$ years!. (See [20], p 119-126 for more details).   \\                   

\noindent                        
$\bullet$ {\it THIRD STEP} \,\, We may finally start with the new operator ${\cal{D}}$ defined by the involutive system $R_4$ with symbol $g_4=0$.The following " {\it fundamental diagram I} " only depends on its left commutative square and $C_0= R_4$. Each horizontal sequence is formally exact and can be constructed step by step. The interest is that we have $C_r= {\wedge}^rT^*\otimes C_0$ because $g_4=0$. It is, even today, not so well known that the three differential sequences appearing in this diagram can be constructed " {\it step by step} " or " {\it as a whole} " ([9],[12],[13],[14]). Accordingly, the reader not familiar with the formal theory of systems of PD equations may find difficult to deal with the following definitions of the Spencer bundles $C_r\subset C_r(E)$ and Janet bundles $F_r$ for an involutive system $R_q \subset J_q(E)$ of order $q$ over $E$: \\
\[  C_r = {\wedge}^rT^* \otimes R_q / \delta ({\wedge}^{r-1}T^*\otimes g_{q+1})  \]
\[  C_r(E) = {\wedge}^rT^*\otimes J_q(E)/ \delta ({\wedge}^{r-1}T^* \otimes S_{q+1}T^*\otimes E)  \]
\[   F_r = {\wedge}^rT^*\otimes J_q(E) / ({\wedge}^rT^*\otimes R_q + \delta ({\wedge}^{r-1}T^*\otimes 
S_{q+1}T^*\otimes E))   \]
For this reason, we prefer to use successive compatibility conditions, starting from the commutative square ${\cal{D}}=\Phi \circ j_4$ on the left of the next diagram. The Janet tabular of the Macaulay system and its prolongations up to order $4$ can be decomposed as follows ([4]):  \\
\[ \left\{ \begin{array}{rccc}
 1  & PDE & order \,\, 4  & class \,\, 3  \\
  4   & PDE & order \,\, 4 & class \,\, 2 \\
 10   &  PDE  & order \,\, 4  & class \,\, 1  \\
   9    &   PDE  & order \,\, 3  &   \\
   3   & PDE & order \,\, 2 &    
\end{array}   \right.
\fbox{  $  \begin{array}{ccc}
1 & 2 & 3 \\
1 & 2 & \bullet  \\
1 & \bullet & \bullet \\
\bullet & \bullet  & \bullet  \\
\bullet & \bullet &  \bullet  
\end{array}  $  }     \]

\noindent
The total number of different single "dots" provides the $4 + 20 + 27 + 9 = 60 $ CC ${\cal{D}}_1$.  \\
The total number of different couples of "dots" provides the $10+27+9=46$ CC ${\cal{D}}_2$.  \\
The total number of different triples of "dots" provides the $ 9 + 3 = 12$ CC ${\cal{D}}_3$.  \\
We obtain therefore the fiber dimensions of the successive Janet bundles in the Janet sequence.\\
The same procedure can be applied to the Spencer bundles in the Spencer sequence by introducing the new $8$ parametric jet indeterminates:  \\ 
\[   z^1=y, z^2=y_1, z^3=y_2 , z^4=y_3, z^5=y_{11}, z^6=y_{12}, z^7=y_{13}, z^8= y_{111}   \]
in the first order system defined by $24$ PD equations ($8$ of class $3$+$8$ of class $2$+$8$ of class $1$):  \\
\[ z^1_1-z^2=0, z^1_2-z^3=0, z^1_3 - z^4=0, ..., z^5_1 - z^8=0, ..., z^6_3- z^8=0, ..., z^7_3=0, ..., 
z^8_3=0 \]

  \[  \begin{array}{rccccccccccccccl}
 &&&&& 0 &&0&&0&  & 0 & & \\
 &&&&& \downarrow && \downarrow && \downarrow  & & \downarrow  &  &\\
  & 0& \longrightarrow& \Theta &\stackrel{j_4}{\longrightarrow}& C_0  &\stackrel{D_1}{\longrightarrow}&  C_1 &\stackrel{D_2}{\longrightarrow} &  C_2 &\stackrel{D_3}{\longrightarrow}  & C_3 & \longrightarrow  & 0  \\
  &&&&& \downarrow & & \downarrow & & \downarrow & & \downarrow &   & \\
   & 0 & \longrightarrow &  E & \stackrel{j_4}{\longrightarrow} & C_0(E)  & \stackrel{D_1}{\longrightarrow} &  C_1(E)  &\stackrel{D_2}{\longrightarrow} & C_2(E) &  \stackrel{D_3}{\longrightarrow} & C_3(E) & \longrightarrow  & 0 \\
   & & & \parallel && \hspace{5mm}\downarrow {\Phi}_0 & &\hspace{5mm} \downarrow {\Phi}_1 & & \hspace{5mm}\downarrow {\Phi}_2 &  &\hspace{5mm}\downarrow {\Phi}_3 & \\
   0 \longrightarrow & \Theta &\longrightarrow &  E & \stackrel{\cal{D}}{\longrightarrow} &  F_0  & \stackrel{{\cal{D}}_1}{\longrightarrow} &  F_1  & \stackrel{{\cal{D}}_2}{\longrightarrow} &  F_2 & \stackrel{{\cal{D}}_3}{\longrightarrow}&  F_3  & \longrightarrow  &  0      \\
   &&&&& \downarrow & & \downarrow & & \downarrow &  & \downarrow &    &   \\
   &&&&& 0 && 0 && 0  &  &  0  &  &
   \end{array}     \]

  \[  \begin{array}{rccccccccccccccl}
 &&&&& 0 &&0&&0&  & 0 & & \\
 &&&&& \downarrow && \downarrow && \downarrow  & & \downarrow  &  &\\
  & 0& \longrightarrow& \Theta &\stackrel{j_4}{\longrightarrow}& 8  &\stackrel{D_1}{\longrightarrow}& 24 &\stackrel{D_2}{\longrightarrow} &  24 &\stackrel{D_3}{\longrightarrow}  & 8 & \longrightarrow  & 0  \\
  &&&&& \downarrow & & \downarrow & & \downarrow & & \downarrow &   & \\
   & 0 & \longrightarrow &  1 & \stackrel{j_4}{\longrightarrow} & 35  & \stackrel{D_1}{\longrightarrow} & 
   84  &\stackrel{D_2}{\longrightarrow} & 70 &  \stackrel{D_3}{\longrightarrow} & 20 & \longrightarrow  & 0 \\
   & & & \parallel && \hspace{5mm}\downarrow {\Phi}_0 & &\hspace{5mm} \downarrow {\Phi}_1 & & \hspace{5mm}\downarrow {\Phi}_2 &  &\hspace{5mm}\downarrow {\Phi}_3 & \\
   0 \longrightarrow & \Theta &\longrightarrow &  1 & \stackrel{\cal{D}}{\longrightarrow} &  27  & \stackrel{{\cal{D}}_1}{\longrightarrow} &  60  & \stackrel{{\cal{D}}_2}{\longrightarrow} &  46 & \stackrel{{\cal{D}}_3}{\longrightarrow}&  12  & \longrightarrow  &  0      \\
   &&&&& \downarrow & & \downarrow & & \downarrow &  & \downarrow &    &   \\
   &&&&& 0 && 0 && 0  &  &  0  &  &
   \end{array}     \]
The morphisms ${\Phi}_1, {\Phi}_2, {\Phi}_3$ in the vertical short exact sequences are inductively induced from the morphism ${\Phi}_0= \Phi$ in the first short exact vertical sequence on the left. The central horizontal sequence can be called " {\it hybrid sequence} " because it is at the same time a Spencer sequence for the first order system $J_5(E) \subset J_1(J_4(E))$ over $J_4(E)$ and a Janet sequence for the involutive injective operator $j_4:E \rightarrow J_4(E)$. It can be constructed step by step, starting with the short exact sequence:   \\
\[  0 \rightarrow J_5(E) \rightarrow J_1(J_4(E)) \rightarrow C_1(E) \rightarrow 0  \]
\[  0 \rightarrow  56   \rightarrow 140  \rightarrow  84  \rightarrow 0  \]
In actual practice, as the system $R_2 \subset J_2(E)$ is homogeneous, it is thus formally integrable and finite type because the system $R_4={\rho}_2(R_2)= ker(\Phi) \subset J_4(E)$ is trivially involutive with a symbol $g_4=0$. Accordingly, ${\cal{D}}=\Phi \circ j_4$ is an involutive operator of order $4$ and we obtain a finite length Janet sequence which is formally exact both on the jet level and on the symbol level, that can only contain the successive first order operators $ {\cal{D}}_1, {\cal{D}}_2, {\cal{D}}_3$. For example, one can determine ${\cal{D}}_2= {\Psi}_2\circ j_1: F_1 \rightarrow F_2$ just by counting the dimensions, either in the long exact jet sequence:    \\
\[  0 \rightarrow R_6 \rightarrow J_6(E) \rightarrow J_2(F_0) \rightarrow J_1(F_1)\stackrel{{\Psi}_2}{ \longrightarrow} F_2 \rightarrow 0  \]
\[   0 \rightarrow 8 \rightarrow 84 \rightarrow 270 \rightarrow 240 \rightarrow dim(F_2) \rightarrow 0 \]
and obtain $dim(F_2)=  - 8 + 84 - 270 + 240=  46  $.   \\
However, one can also use the fact that $dim(E)=1$ and $g_4=0 \Rightarrow g_6=0$ while introducing the restriction $\sigma({\Psi}_2)$ of ${\Psi}_2$ to $T^*\otimes F_1 \subset J_1(F_1)$ in the long exact symbol sequence:  \\
\[   0 \rightarrow S_6T^* \rightarrow S_2T^*\otimes F_0 \rightarrow T^*Ê\otimes F_1 \stackrel{\sigma({\Psi}_2)}{\longrightarrow} F_2 \rightarrow 0  \]
\[   0 \rightarrow 28 \rightarrow 162 \rightarrow 180 \rightarrow dim(F_2) \rightarrow 0  \]
in order to obtain again $ dim(F_2)= 28 - 162 + 180 =  46$.  \\

We wish good luck to anybody using Computer Algebra because one should have to deal with a matrix $540 \times 600$ in order to describe the prolongation morphism $ J_3(F_0) \rightarrow J_2(F_1)$. Nevertheless, in order to give a hint, we recall the vanishing of the Euler-Poincar\'{e} characteristic as we can check successively:   \\
\[       8 - 24 + 24 - 8 = 0, \,\,\, -1 + 35 - 84 + 70 - 20 = 0, \,\,\, -1 + 27 - 60 + 46 - 12 =0  \]
In the case of finite type systems, the usefulness of the Spencer sequence is so evident, like on such an example, that it needs no comment.  \\
We invite the reader to treat separately but similarly the system:  \\
 \[  y_{33} - y_{11}=0,\,\,  y_{23}=0, \,\, y_{22} - y_{11}=0  \]
and to compare the various extension modules.  \\

\noindent
{\bf 3) SOLUTION}   \\
According to the previous sections, it only remains to consider the two cases $n=1$ and $n=2$. For simplicity, we shall only consider the situation of the Euclidean metric and the corresponding linear systems. We let the reader treat by himself the nonlinear counterparts.    \\

\noindent
$ \bullet  \,\, {\it CASE\,\,\,\,  n=1}$ \\
With $\omega \neq 0$, we may consider a section ${\xi}_3= (\xi(x), {\xi}_x(x),{\xi}_{xx}(x),{\xi}_{xxx})$ and introduce the classical Killing system $R_1\subset J_1(T)$  by means of the {\it formal Lie derivative}:  \\
\[     \Omega \equiv L({\xi}_1)\omega \equiv 2 \omega {\xi}_x + \xi {\partial}_x\omega=0  \]
Similarly, with the Christoffel symbol $\gamma = \frac{1}{2\omega} {\partial}_x\omega$, we may consider:  \\
\[   \Gamma \equiv L({\xi}_2)\gamma \equiv {\xi}_{xx} + \gamma {\xi}_x + \xi {\partial}_x \gamma =0   \]
The conformal Killing system can be defined with a conformal factor as:  \\
\[   \Omega \equiv L({\xi}_1)\omega \equiv 2 \omega {\xi}_x + \xi {\partial}_x\omega=2 A(x) \omega  \]
and its first prolongation becomes: \\
\[  \Gamma \equiv L({\xi}_2)\gamma \equiv {\xi}_{xx} + \gamma {\xi}_x + \xi {\partial}_x \gamma = A_x(x) \]
The elimination of $A(x)$ or $A_x(x)$ does not provide any OD equation of order $1$ or $2$. Moreover, we let the reader check that ${\xi}_2=j_2(\xi) \Rightarrow {\partial}_xA(x) - A_x(x)=0$ as a way to understand the part plaid by the Spencer operator and the reason for introducing $2A(x)$. With more details, dividing the Killing system by $2\omega$, we get ${\xi}_x + \gamma \xi= A(x)$. Differentiating this  OD equation, we get:  \\
\[    {\partial}_x{\xi}_x + \gamma {\partial}_x\xi + {\partial}_{x}\gamma  \xi =  {\partial}_x A(x)   \]
and we just need to substract the OD equation $\Gamma = A_x(x)$ in order to get:  \\
\[         ({\partial}_x{\xi}_x - {\xi}_{xx}) + \gamma ({\partial}_x\xi - {\xi}_x) = {\partial}_xA(x) - A_x(x)  \]
In order to escape from the previous situation while having a vanishing symbol $g_3=0$, we may consider the new {\it unusual prolongation}: \\
\[   {\xi}_{xxx} + \gamma {\xi}_{xx} + 2 ({\partial}_x\gamma) {\xi}_x + \xi {\partial}_{xx} \gamma = 0  \]
and substract the second order OD equation $\Gamma = 0$ multiplied by $\gamma$ while introducing the new geometric object $\nu={\partial}_x\gamma - \frac{1}{2}{\gamma}^2$ in order to obtain the third order infinitesimal Lie equation: \\
\[    L({\xi}_3)\nu \equiv  {\xi}_{xxx} + 2 \nu {\xi}_x + \xi {\partial}_x \nu = 0  \]
The nonlinear framework, not known today because the work of Vessiot is still not acknowledged, explains the successive inclusions $\gamma \in j_1(\omega), \nu \in j_1(\gamma)$. Indeed, if we consider the {\it translation group} $(y=x + a, a=cst)$ and the bigger {\it isometry group} $(y=x + a, y = - x + a, a=cst)$, the inclusion of groups of the real line:  \\
\[    translation\,\, group \subset isometry \,\, group \subset affine \,\,  group  \subset projective \,\, group   \]
with the respective finite Lie equations in Lie form with the jet coordinates $(x,y, y_x, y_{xx}, y_{xxx})$:   \\
\[  {\alpha}(y)y_x={\alpha}(x), \,\,\,\,\,   {\omega}(y) (y_x)^2 = {\omega}(x) , \,\,\,\,\,      \frac{y_{xx}}{y_x} + \gamma (y) y_x = \gamma (x),                              
     \,\,\,\,\,   \frac{y_{xxx}}{y_{x}} - \frac{3}{2} (\frac{y_{xx}}{y_x})^2 + \nu (y) (y_x)^2 = \nu (x)   \]
where we recognize the Schwarzian third order differential invariant of the {\it projective group}.   \\
Of course, we have $\alpha = 1 \Rightarrow \omega =1 \Rightarrow \gamma =0  \Rightarrow  \nu =0  $ and the respective linearizations:  \\
\[   y_x=1 \,\, \Rightarrow \,\, {\xi}_x=0, \,\,\, y_{xx}=0 \,\, \Rightarrow \,\, {\xi}_{xx} = 0 , \,\,\, 
       \frac{y_{xxx}}{y_x} - \frac{3}{2} (\frac{y_{xx}}{y_x})^2=0 \,\, \Rightarrow \,\, {\xi}_{xxx}=0   \]

The Janet tabular of the conformal system order $3$ can be decomposed as follows:  \\
\[ \left\{ \begin{array}{rccc}
1 & PDE & order \,\, 3  & class \,\, 1  
\end{array}   \right.
\fbox{  $  \begin{array}{c}
1      
\end{array}  $  }     \]

\noindent
The total number of different single "dots" provides the $0 $ CC ${\cal{D}}_1$.  \\
We obtain therefore the fiber dimensions of the successive Janet bundles in the Janet sequence.\\
The same procedure can be applied to the other canonical differential sequences.  \\ 
When $n=1$, one has $3$ parameters ($1$ translation + $1$ dilatation + $1$ elation) and we get the following " {\it fundamental diagram I} "  only depending on the left commutative square:  \\

\[  \begin{array}{rccccccccc}
  & && & &0 & & 0 &  &  \\
   & & & & & \downarrow & & \downarrow & &  \\
   & 0 &\longrightarrow &\Theta &\stackrel{j_3}{\longrightarrow} & 3 &\stackrel{D_1}{\longrightarrow} &3 &\longrightarrow 0 & \hspace{5mm} Spencer   \\
   & & & & & \downarrow & & \parallel  &  &  \\
    & 0 & \longrightarrow &1 &\stackrel{j_3}{\longrightarrow} & 4 &\stackrel{D_1}{\longrightarrow} & 3 & \longrightarrow 0&   \\
    & & & \parallel & & \hspace{3mm}\downarrow \Phi & & \downarrow & &  \\
    0\longrightarrow & \Theta & \longrightarrow & 1 & \stackrel{{\cal{D}}}{\longrightarrow} & 1 & \longrightarrow & 0 & & \hspace{5mm} Janet  \\
    & & & & & \downarrow & & & &  \\
     & & & & & 0 & & & & 
 \end{array}  \]

\noindent 
In this diagram, the operator $j_3:\xi(x)\rightarrow (\xi(x)=\xi(x),{\partial}_x\xi(x)={\xi}_x(x),{\partial}_{xx}\xi(x)={\xi}_{xx}(x), {\partial}_{xxx}\xi(x)={\xi}_{xxx}(x))$ has compatibility conditions $D_1{\xi}_3=0$ induced by $d$ and the space of solutions $\Theta $ of ${\cal{D}}=\Phi\circ j_3: \xi(x) \rightarrow {\partial}_{xxx} \xi (x)$ is generated over the constants by the three  infinitesimal generators:  \\
\hspace*{2cm}${\theta}_1={\partial}_x$ (translation), \,\,\, ${\theta}_2= x{\partial}_x$ (dilatation), \,\,\,  ${\theta}_3=\frac{1}{2}x^2{\partial}_x$ (elation) \\  
 of the action and coincides with the projective group of the real line. \\

\noindent
$  \bullet \,\, {\it CASE \,\,\,\, n=2}$ \\
The classical approach is to consider the infinitesimal conformal Killing system for $n=2$ and eliminate the infinitesimal conformal factor $2A(x)$ as follows by introducing the {\it formal} and the {\it effective} Lie derivatives such that $L(j_1(\xi))= {\cal{L}}(\xi)$:   \\
\[ \Omega \equiv L({\xi}_1)\omega =2A(x)\omega \Rightarrow 
{\xi}^1_1=A(x), {\xi}^1_2 + {\xi}^2_1=0, {\xi}^2_2=A(x) \Rightarrow 
{\xi}^2_2 - {\xi}^1_1=0, {\xi}^1_2 +{\xi}^2_1=0   \]
that is to say the elimination of $A$ is just producing locally the two well known Cauchy-Riemann equations allowing to define infinitesimal complex transformations of the plane, that is to say an infinite dimensional Lie pseudogroup which is by no way providing a finite dimensional Lie group. As such an operator has no {\it compatibility condition} (CC), we obtain by one prolongation $2\times 2=4$ second order equations but another prolongation does not provide a zero symbol at order $3$ and it is just such a delicate step that we have to overcome by adding $2\times 4=8$ homogeneous third order PD equations. The {\it only possibility} which is coherent with homogeneity is thus to consider the following system and to prove that it is defining a  system of infinitesimal Lie equations leading to $2\times (1+2+3+4) - (2+4+8)=20-14=6$ infinitesimal generators.  \\
\[  \left \{ \begin{array}{c}
  {\xi}^k_{ijr}=0   \\
   {\xi}^2_{22} - {\xi}^1_{12}=0, {\xi}^1_{22} + {\xi}^2_{12}=0, {\xi}^2_{12} - {\xi}^1_{11}=0, 
{\xi}^1_{12} + {\xi}^2_{11}=0   \\
  {\xi}^2_2 - {\xi}^1_1=0, {\xi}^1_2 + {\xi}^2_1=0   
  \end{array} \right.   \]
where the $4$ second order PD equations can also be rewritten with $\Delta=d_{11} + d_{22}$ as:  \\
\[  \Delta {\xi}^2\equiv {\xi}^2_{22} +{\xi}^2_{11}=0, \Delta {\xi}^1 \equiv {\xi}^1_{22} + {\xi}^1_{11}=0,
     {\xi}^2_{12} - {\xi}^1_{11}=0, {\xi}^1_{12} + {\xi}^2_{11}=0  \]
     
The general solution of the $8$ third order PD equations can be written with $12$ arbitrary constant parameters as:   \\ 
\[  {\xi}^1= \frac{1}{2}a(x^1)^2 + bx^1x^2 + \frac{1}{2}c(x^2)^2 + dx^1 + ex^2 + f   \]   
\[  {\xi}^2= \frac{1}{2}\bar{a}(x^1)^2 +\bar{b}x^1x^2 + \frac{1}{2}\bar{c}(x^2)^2 + \bar{d}x^1 + 
\bar{e}x^2 + g  \]
Taking into account the first and second order PD equations, we must have the relations:  \\
\[    \bar{b}=a, \bar{c}=b, \bar{a} + b=0, \bar{b} + c= 0, \bar{e} = d, \bar{d} + e=0  \] 
and the final number of parameters is indeed reduced to $2 + 1 + 1 + 2=6$ arbitrary parameters. Collecting the above results, we obtain the 6 infinitesimal generators:  \\
\[ a \rightarrow \frac{1}{2}((x^1)^2  - (x^2)^2){\partial}_1 + x^1x^2{\partial}_2  \]
\[ b \rightarrow x^1x^2 {\partial}_1+ \frac{1}{2}( (x^2)^2 - (x^1)^2){\partial}_2  \]
\[ - e \rightarrow x^1{\partial}_2 - x^2{\partial}_1, \,\, d\rightarrow x^1 {\partial}_1 + x^2{\partial}_2 \]
 \[   f \rightarrow {\partial}_1, \,\, g \rightarrow {\partial}_2   \]
We find back the two infinitesimal generators of the elations, namely:   \\
\[    {\theta}_1= - \frac{1}{2}((x^1)^2 + (x^2)^2){\partial}_1 + x^1(x^1{\partial}_1 +x^2 {\partial}_2) =
\frac{1}{2}((x^1)^2 - (x^2)^2){\partial}_1 + x^1x^2{\partial}_2  \]
and ${\theta}_2$ obtained by exchanging $x^1$ with $x^2$. \\  
Contrary to the situation met when $n\geq 3$ where one starts with a groupoid of order $1$ and obtains groupoids of order $2$ or $3$ after one or two prolongations, in the present situation, we have to check directly the commutation relations for the six infinitesimal generators already found, namely:  \\
\[      [{\partial}_1, {\theta}_1] = x^1{\partial}_1 + x^2 {\partial}_2, \,\,\,
   [{\partial}_2, {\theta}_1] = x^1{\partial}_2 - x^2{\partial}_1\]
\[   [x^1{\partial}_2 - x^2{\partial}_1,{\theta}_1] = - {\theta}_2, \,\,\, 
   [x^1{\partial}_1 + x^2{\partial}_2, {\theta}_1] = {\theta}_1, \,\,\,
[{\theta}_1 , {\theta}_2 ] = 0  \]

We have thus obtained in an unexpected way the desired $2$ translations, $1$ rotation, $1$ dilatation and $2$ elations of the conformal group when $n=2$.     \\

At order one, we may consider the classical Killing system $R_1$ obtained by preserving $\omega$, the Weyl system ${\tilde{R}}_1$ and the conformal system ${\hat{R}}_1$ with $R_1 \subset {\tilde{R}}_1={\hat{R}}_1\subset J_1(T)$ and $dim({\tilde{R}}_1/R_1)=1$. At order two, we have the strict inclusions $R_2 \subset {\tilde{R}}_2 \subset {\hat{R}}_2$ with $R_2={\rho}_1(R_1)$ preserving 
$(\omega,\gamma)\simeq j_1(\omega)$, ${\tilde{R}}_2 \subset {\rho}_1({\tilde{R}}_1)$ obtained by preserving $(\hat{\omega}, \gamma)$ and ${\hat{R}}_2={\rho}_1({\hat{R}}_1)$ obtained by preserving 
$ (\hat{\omega},\hat{\gamma})\simeq j_1(\hat{\omega})$. The main difference with the case $n\geq 3$ is that {\it now} $R_3= {\rho}_2(R_1)$ has a symbol $g_3=0$, ${\tilde{R}}_3={\rho}_1({\tilde{R}}_2)$ has also a symbol ${\tilde{g}}_3=0$ {\it but} that ${\hat{R}}_3 \subset {\rho}_1({\hat{R}}_2)$ with strict inclusion in order to have {\it now} ${\hat{g}}_3=0$, even though ${\rho}_1({\hat{g}}_2)\neq 0$. However, we are now able to deal with three {\it trivially involutive} systems having zero symbols and we have the strict inclusions $R_3 \subset {\tilde{R}}_3 \subset {\hat{R}}_3$ with respective dimensions $3 < 4 < 6$ according to the basic inequalities $n(n+1)/2 < (n(n+1)/2)+1 < (n+1)(n+2)/2$ valid in arbitrary dimension $n\geq 1$. The interest of this result is that we have for the Spencer bundles the strict inclusions 
$C_0 \subset {\tilde{C}}_0 \subset {\hat{C}}_0$ of the zero Spencer bundles, leading to the strict inclusions of the respective linear Spencer sequences because:   \\
\[  g_3 ={\tilde{g}}_3={\hat{g}}_3=0 \Rightarrow  C_r={\wedge}^rT^*\otimes C_0, {\tilde{C}}_r={\wedge}^rT^*\otimes {\tilde{C}}_0, {\hat{C}}_r= {\wedge}^rT^* \otimes {\hat{C}}_0 \Rightarrow  C_r \subset {\tilde{C}}_r \subset {\hat{C}}_r  \] 
in agrement with many recent results ([21],[22],[23],[24]). As in Example $2.2$, we let the reader introduce the $6$ parametric jet indeterminates $z^1=y^1,z^2=y^2,z^3=y^1_1,z^4=y^2_1,z^5=y^1_{11},z^6=y^2_{11}$.  \\

The Janet tabular of the conformal Killing system and its prolongations up to order $3$ can be decomposed as follows:  \\
\[ \left\{ \begin{array}{rccc}
 2 & PDE & order \,\, 3  & class \,\, 2  \\
 6 & PDE & order \,\, 3 & class \,\, 1 \\
 4     &  PDE  & order \,\, 2  &         \\
 2       &   PDE  & order \,\, 1  &    
\end{array}   \right.
\fbox{  $  \begin{array}{cc}
1 & 2   \\
1 & \bullet  \\
\bullet & \bullet  \\
\bullet & \bullet 
\end{array}  $  }     \]

\noindent
The total number of different single "dots" provides the $6 + 8 + 4 = 18 $ CC ${\cal{D}}_1$.  \\
The total nuber of different couples of "dots" provides the $4 + 2 = 6$ CC ${\cal{D}}_2$.  \\
We obtain therefore the fiber dimensions of the successive Janet bundles in the Janet sequence.\\
The same procedure can be applied to the other canonical differential sequences.  \\ 
When $n=2$, one has $6$ parameters ($2$ translations + $1$ rotation + $1$ dilatation + $2$ elations) and we get the following " {\it fundamental diagram I} "  only depending on the left commutative square: \\
  
  \[  \begin{array}{rccccccccccccc}
 &&&&& 0 &&0&&0&  & \\
 &&&&& \downarrow && \downarrow && \downarrow  &\\
  & 0& \longrightarrow& \Theta &\stackrel{j_3}{\longrightarrow}& 6  &\stackrel{D_1}{\longrightarrow}&  12 &\stackrel{D_2}{\longrightarrow} &  6 &\longrightarrow  0 & \hspace{3mm}Spencer  \\
  &&&&& \downarrow & & \downarrow & & \downarrow & &    & \\
   & 0 & \longrightarrow &  2 & \stackrel{j_3}{\longrightarrow} & 20  & \stackrel{D_1}{\longrightarrow} &  
   30  &\stackrel{D_2}{\longrightarrow} & 12 &   \longrightarrow 0 &\\
   & & & \parallel && \hspace{5mm}\downarrow {\Phi}_0 & &\hspace{5mm} \downarrow {\Phi}_1 & & \hspace{5mm}\downarrow {\Phi}_2 &  &\\
   0 \longrightarrow & \Theta &\longrightarrow &  2 & \stackrel{\cal{D}}{\longrightarrow} &  14  & \stackrel{{\cal{D}}_1}{\longrightarrow} &  18  & \stackrel{{\cal{D}}_2}{\longrightarrow} &  6 & \longrightarrow  0 & \hspace{5mm} Janet       \\
   &&&&& \downarrow & & \downarrow & & \downarrow &      &\\
   &&&&& 0 && 0 && 0  &  &
   \end{array}     \]

\noindent
$ \bullet  \,\, {\it CASE \,\,\,\, n=3}$  \\
The Janet tabular of the conformal Killing system and its prolongations up to order $3$ can be decomposed as follows:  \\
\[ \left\{ \begin{array}{rccc}
 3 & PDE & order \,\, 3  &  class \,\, 3 \\
 9 & PDE & order \,\, 3  & class \,\, 2  \\
18 & PDE & order \,\, 3 & class \,\, 1 \\
15     &  PDE  & order \,\, 2  &         \\
 5       &   PDE  & order \,\, 1  &    
\end{array}   \right.
\fbox{  $  \begin{array}{ccc}
1 & 2 & 3  \\
1 & 2 & \bullet  \\
1 & \bullet & \bullet  \\
\bullet & \bullet  & \bullet  \\
\bullet & \bullet & \bullet
\end{array}  $  }     \]

\noindent
The total number of different single "dots" provides the $  9 + 36 + 45 + 15 = 105 $ CC ${\cal{D}}_1$.  \\
The total number of different couples of "dots" provides the $ 18 + 45 + 15 = 78 $ CC ${\cal{D}}_2$.  \\
The total number of different triples of "dots"  provides the $15 + 5 =20 $ CC ${\cal{D}}_3$. \\
We obtain therefore the fiber dimensions of the successive Janet bundles in the Janet sequence.  \\
The same procedure can be applied to the other canonical differential sequences and we get the desired " {\it fundamental diagram I} "  below:   \\ 

  \[  \begin{array}{rccccccccccccccc}
 &&&&& 0 &&0&&0&  & 0 \\
 &&&&& \downarrow && \downarrow && \downarrow  & &\downarrow  &   \\
  & 0& \longrightarrow& \Theta &\stackrel{j_3}{\longrightarrow}& 10  &\stackrel{D_1}{\longrightarrow}&  30 &\stackrel{D_2}{\longrightarrow} &  30 & \stackrel{D_3}{\longrightarrow}  &10 &  \longrightarrow 0 \\
  &&&&& \downarrow & & \downarrow & & \downarrow & & \downarrow   & &     \\
   & 0 & \longrightarrow &  3 & \stackrel{j_3}{\longrightarrow} & 60  & \stackrel{D_1}{\longrightarrow} &  
   135  &\stackrel{D_2}{\longrightarrow} & 108 & \stackrel{D_3}{\longrightarrow} & 30 & \longrightarrow 0 \\
   & & & \parallel && \hspace{5mm}\downarrow {\Phi}_0 & &\hspace{5mm} \downarrow {\Phi}_1 & & \hspace{5mm}\downarrow {\Phi}_2 &  & \hspace{5mm} \downarrow {\Phi}_3 &  &     \\
   0 \longrightarrow & \Theta &\longrightarrow &  3 & \stackrel{\cal{D}}{\longrightarrow} &  50  & \stackrel{{\cal{D}}_1}{\longrightarrow} &  105  & \stackrel{{\cal{D}}_2}{\longrightarrow} &  78     &      
   \stackrel{{\cal{D}}_3}{\longrightarrow}&  20   & \longrightarrow 0       \\
   &&&&& \downarrow & & \downarrow & & \downarrow &    & \downarrow  &    &      \\
   &&&&& 0 && 0 && 0  &  & 0  &   &
   \end{array}     \]
   
   \noindent
We have $10$ parameters ( $3$ translations, $3$ rotations, $1$ dilataion, $3$ elations).  \\
The computation of $dim(C_3(E))=30$ needs to determine the rank of a $1200 \times 1350$ matrix !.\\ \\

\noindent
{\bf 4) CONCLUSION}  \\  

We have shown that the true important specific property of the conformal group, at least for applications to physics, is that, even if it is defined as a specific Lie pseudogroup of transformations, it is in fact a Lie group of transformations with a finite number $(n+1)(n+2)/2$ of parameters or infinitesimal generators wen $n\geq 3$. Accordingly, in dimension $n=1$, we have no OD equation of order $1$ {\it and} $2$, a result leading therefore to add $ 1$ unexpected OD equation of order $3$. Similarly, when $n=2$, we obtain the  Cauchy-Riemann PD equations defining an infinite dimensional Lie pseudogroup and we have therefore to add, again in a totally unexpected way, as many third order PD equations as the number of jet coordinates of strict order $3$. When $n=3$, the fact that the analogue of the Weyl operator for describing the CC of the conformal operator is of order $3$ is rather unpleasent but this is nothing compared to the fact that, when $n=4$, the analogue of the Bianchi operator for describing the CC of the  previous second order CC playing the part of the Weyl CC is of order $2$ {\it again}. And we don't speak about the case $n=5$ ([18],[20]). Though these results can be checked by means of computer algebra and are confirmed by the use of the fundamental diagram I, they do not seem to be known today. Accordingly, any {\it physical theory} (existence of gravitational waves or black holes, ... ) which is not coherent with {\it differential homological algebra} (vanishing of the first and second extension modules for the Poincar\'{e} sequence in the previous examples, ...) must be revisited in the light of these new {\it mathematical tools}, even if it seems apparently well established ([16],[21],[22]). \\ \\

\noindent
{\bf REFERENCES}  \\

\noindent
[1] Bj\"{o}rk, J.-E.: Analytic ${\cal{D}}$-Modules and Applications, Mathematics and Its Applications, 247, Kluwer, 1993.   \\
\noindent
[2] Eisenhart, L.P.: Riemannian Geometry, Princeton University Press, Princeton, 1926.\\
\noindent
[3] Gasqui, J., Goldschmidt, H.: D\'{e}formations Infinit\'{e}simales des Structures Conformes Plates, Progress in Mathematics, Vol. 52, Birkhauser, Boston, 1984.   \\
\noindent
[4] Janet, M.: Sur les Syst\`{e}mes aux D\'{e}riv\'{e}es Partielles, Journal de Math., 8 (1920) 65-151. \\
\noindent 
[5] Kashiwara, M.: Algebraic Study of Systems of Partial Differential Equations, M\'emoires de la Soci\'et\'e Math\'ematique de France, 63, 1995) (Transl. from Japanese of his 1970 Master's Thesis).\\
\noindent
[6] Kumpera, A., Spencer, D.C.: Lie Equations, Ann. Math. Studies 73, Princeton University Press, Princeton, 1972.\\
\noindent
[7] Macaulay, F.S.: The Algebraic Theory of Modular Systems, Cambridge Tracts, vol. 19, Cambridge University Press, London, 1916. Stechert-Hafner Service Agency, New-York, 1964.\\
\noindent
[8] Northcott, D.G.: An Introduction to Homological Algebra, Cambridge university Press, 1966.  \\
\noindent
[9] Pommaret, J.-F.: Systems of Partial Differential Equations and Lie Pseudogroups, Gordon and Breach, New York (1978); Russian translation: MIR, Moscow, 1983.\\
\noindent
[10] Pommaret, J.-F.: Differential Galois Theory, Gordon and Breach, New York, 1983.\\
\noindent
[11] Pommaret, J.-F.: Lie Pseudogroups and Mechanics, Gordon and Breach, New York, 1988.\\
\noindent
[12] Pommaret, J.-F.: Partial Differential Equations and Group Theory, Kluwer, 1994.\\
https://doi.org/10.1007/978-94-017-2539-2    \\
\noindent
[13] Pommaret, J.-F.: Partial Differential Control Theory, Kluwer, Dordrecht, 2001.\\
\noindent
[14] Pommaret, J.-F.: Algebraic Analysis of Control Systems Defined by Partial Differential Equations, in "Advanced Topics in Control Systems Theory", Springer, Lecture Notes in Control and Information Sciences 311, 2005, Chapter 5, pp. 155-223.\\
\noindent
[15] Pommaret, J.-F.: Parametrization of Cosserat Equations, Acta Mechanica, 215 (2010) 43-55.\\
https://doi.org/10.1007/s00707-010-0292-y  \\
\noindent
[16] Pommaret, J.-F.: The Mathematical Foundations of General Relativity Revisited, Journal of Modern Physics, 4 (2013) 223-239. \\
 https://doi.org/10.4236/jmp.2013.48A022   \\
 \noindent
[17] Pommaret, J.-F.: Relative Parametrization of Linear Multidimensional Systems, Multidim. Syst. Sign. Process., 26 (2015) 405-437.  \\
DOI 10.1007/s11045-013-0265-0   \\
\noindent
[18] Pommaret, J.-F.: Airy, Beltrami, Maxwell, Einstein and Lanczos Potentials revisited, Journal of Modern Physics, 7 (2016) 699-728. \\
\noindent
[19] Pommaret,J.-F.:From Thermodynamics to Gauge Theory: the Virial Theorem Revisited, pp. 1-46 in "Gauge Theories and Differential geometry,", NOVA Science Publisher (2015).  \\
\noindent
https://doi.org/10.4236/jmp.2016.77068   \\
\noindent
[20] Pommaret, J.-F.: Deformation Theory of Algebraic and Geometric Structures, Lambert Academic Publisher (LAP), Saarbrucken, Germany, 2016. A short summary can be found in "Topics in Invariant Theory ", S\'{e}minaire P. Dubreil/M.-P. Malliavin, Springer Lecture Notes in Mathematics, 1478 (1990) 244-254.\\
https://arxiv.org/abs/1207.1964  \\
\noindent
[21] Pommaret, J.-F.: Why Gravitational Waves Cannot Exist, J. of Modern Physics, 8 (2017) 2122-2158.  \\
https://doi.org/104236/jmp.2017.813130    \\
\noindent
[22] Pommaret, J.-F.: New Mathematical Methods for Physics, Mathematical Physics Books, Nova Science Publishers, New York, 2018, 150 pp. \\
\noindent
[23] Pommaret, J.-F.: Generating Compatibility Conditions and General Relativity, J. of Modern Physics, 10, 3 (2019) 371-401.  \\
https://doi.org/10.4236/jmp.2019.103025   \\
\noindent
[24] Pommaret, J.-F.: Differential Homological Algebra and General Relativity, J. of Modern Physics, 10 (2019) 1454-1486. \\
https://doi.org/10.4236/jmp.2019.1012097   \\
\noindent
[25] Rotman, J.J.: An Introduction to Homological Algebra, Academic Press, 1979.   \\
\noindent
[26] Schneiders, J.-P.: An Introduction to ${\cal{D}}$-Modules, Bull. Soc. Roy. Sci. Li\'{e}ge, 63 (1994) 223-295.  \\
\noindent
[27] Spencer, D.C.: Overdetermined Systems of Partial Differential Equations, Bull. Am. Math. Soc., 75 (1965) 1-114.\\
\noindent
[28] Vessiot, E.: Sur la Th\'{e}orie des Groupes Infinis, Ann. Ec. Norm. Sup., 20 (1903) 411-451.   \\ 
\noindent

\end{document}